\documentstyle[11pt,amssymb]{amsart}
\def\limfunc#1{\hbox{\rm #1}}

\newtheorem{theorem}{Theorem}
\newtheorem{lemma}[theorem]{Lemma}

\newtheorem{corollary}[theorem]{Corollary}

\begin{document}

\title[Torsion modules, lattices and p-points]
      {Torsion modules, lattices and p-points}

\author{Paul C. Eklof}
\address{Department of Mathematics, University of California, Irvine;
Irvine, CA 92697, USA}                                              

\author{Birge Huisgen-Zimmermann}
\address{Department of Mathematics, University of California, Santa
Barbara, Santa Barbara, CA 93106, USA} 

\author{Saharon Shelah}
\address{Mathematics Institute, Hebrew University, Jerusalem 91904, Israel}

\thanks{Second author  partly supported by NSF grant DMS
95-00453. Third author supported by German-Israeli Foundation for
Scientific Research \& Development Grant No. G-294.081.06/93. Pub. No. 617}

\maketitle

%% Note that there is no paragraph spacing between two affiliations
%% on one line

%start of body        

\begin{abstract} Answering a long-standing question in the
theory of  torsion modules, we show that weakly productively
bounded domains are  necessarily productively bounded. (See
the introduction for  definitions.)  Moreover, we prove a
twin result for  the ideal lattice  $L$  of a domain equating
weak and strong  global intersection conditions for families
$(X_i)_{i
 \in  I}$ of subsets of $L$ with the  property that
$\bigcap_{i  \in  I} A_i
\ne 0$ whenever $A_i  \in  X_i$.   Finally, we show that, for
domains with Krull dimension (and countably  generated
extensions thereof), these lattice-theoretic conditions are
equivalent to  productive boundedness.
\end{abstract}

\section{Introduction}

This paper continues a series of articles (see \cite{gz},
\cite{bg},\cite{bz},\cite{z2}) dealing with the following
problem about torsion modules:
\begin{quote} For which (right) Ore domains $R$ is  it true
that, given any family $(M_i)_{i  \in  I}$ of torsion (right)
$R$-modules, the intersection of annihilators, $\bigcap_{i
\in  J}
\limfunc{ann}(M_i)$, is nonzero for some cofinite subset
$J \subseteq I$ provided that the direct product $\prod _{i
\in  I}M_i$ is  torsion?
\end{quote}
\noindent   In other words, over which Ore domains is the
obvious sufficient condition for a direct product of torsion
modules to be torsion also necessary?  We call such Ore
domains (right) {\it productively bounded}.  The fact that
Dedekind domains are productively bounded plays a pivotal
role in the theory of direct-sum decompositions of direct
products of modules over such domains \cite{z1}, which
triggered interest in the  property for more general rings.
The known classes of productively  bounded domains include
all Ore domains of countable Krull dimension \cite{gz}  and
all commutative noetherian domains \cite{bg}.

As a consequence of our main  results, these classes of
positive examples can be enlarged: namely each  commutative
domain which is a countably generated extension of a
productively bounded domain inherits this property (Section
4).

Already in the first investigation of the topic (\cite{gz}),
it turned out that the following weakened condition on $R$,
accordingly labeled (right) {\it weak productive
boundedness}, is far more accessible in concrete situations:
whenever $(M_i)_{i  \in  I}$ is a family of $R$-modules such
that $\prod_{i  \in  I} M_i$ is torsion, the annihilators
$\limfunc{ann}_R(M_i)$ are nonzero for all but finitely many
$i  \in  I$. This naturally raised the question whether weak
productive boundedness implies the full boundedness condition
in general.  Our main theorem answers this question
affirmatively (Section 1).  The proof rests on the notion of
a p-point, a certain type of ultrafilter on $\omega$ (defined
below).

There is an immediately neighboring pair of properties for
the dual ideal  lattice of $R$ which, at least on the face of
it, is somewhat stronger than the  two boundedness conditions
for $R$ discussed so far.  Consider any lattice
$L$ which is complete and finitely join-irreducible, the
latter meaning  that its largest element 1 is not the join of
two strictly smaller  elements. (Observe that the dual  ideal
lattice of an Ore domain  satisfies these conditions.)   If
$(X_i)_{i \in  I}$ is a family of nonempty subsets of $L$, an
element $x= (x_i)_{i \in  I}$ of $\prod_{i \in  I}X_i$ is
referred to as a {\it transversal  of $  (X_i)_{i \in  I}$ }
and such a transversal is said to be {\it bounded} if
$\bigvee \{x_i:i \in  I\}<1$. We say that $L$ is {\it
uniformly transversally bounded }if for  each family
$(X_i)_{i \in  I}$ of nonempty subsets of $L$ all of whose
transversals are bounded, there is a cofinite subset $J$ of
$I$ such that $\bigvee (\bigcup_{i \in  J}X_i)<1$; moreover,
we call $L$ {\it
 tranversally bounded} if for each family $(X_i)_{i \in  I}$
of nonempty subsets of $L$, all of whose transversals are
bounded, $\bigvee  X_i<1$ for all but finitely many $i  \in
I$.

Given a domain $R$, we denote by $L_R$ the dual of its right
ideal  lattice.  It is easy to see that transversal
boundedness of $L_R$ (resp.  uniform transversal boundedness
of $L_R$) entails weak productive  boundedness (resp.
productive boundedness) of $R$.  It remains open  whether the
reverse implications hold in general.  However, we show in
Section 3 that for domains with Krull dimension in the sense
of Gordon  and Robson \cite{gr}  all of these boundedness
conditions  are  equivalent. This completes a round-trip from
torsion modules to  ultrafilters through lattices which began
with the papers we listed at  the outset.

For the dual ideal lattices of {\it arbitrary} domains we
prove transversal  boundedness to be equivalent to uniform
transversal boundedness (Section  2).  Interestingly, this
equivalence distinguishes  dual ideal lattices  of domains
from abstract lattices.  In fact, the second author showed
that the continuum hypothesis guarantees the existence of
lattices which  are transversally bounded, but not uniformly
so (see \cite{z2}).  In Section 2 of the present paper, we
observe that this conclusion is actually independent of ZFC,
even independent of ZFC plus the negation of the continuum
hypothesis.

Throughout, we shall assume that $R$ is a commutative
integral  domain. However, all of  our results in Sections
1--3 actually carry over to arbitrary Ore domains  in a quite
obvious manner (for module read ``right module" and for ideal
read ``right ideal"). This is not true for the theorems of
Section  4, where we warn the reader about this point.

 That $R$ is productively bounded (resp. weakly productively
bounded) will be abbreviated by `$R$ is PB' (resp. `$R$ is
wPB').   Moreover, transversal boundedness and uniform
transversal boundedness  will be denoted by `TB' and `UTB',
respectively, for short.

We recall the definitions of some special kinds of
ultrafilters on
$\omega$ (see \cite[pp 257-9]{J}, for example).  An
ultrafilter $U$ on
$\omega$ is called a {\it p-point} (resp. a {\it Ramsey}
ultrafilter) if,  whenever $\omega$ is expressed as a
disjoint union $\coprod_{k  \in
\omega} N_k$ of  subsets $N_k$, none of which belongs to $U$,
there is a  set $Y  \in  U$ such that $Y \cap N_k$ is finite
(resp. $|Y \cap N_k| \le
 1$) for all $k  \in  \omega$.  It is known that either of CH
or (MA +
 $\neg$ CH) implies the existence of Ramsey ultrafilters.  On
the other hand, by a result due to the third author
\cite{Sh}, the non-existence of p-points is consistent with
ZFC + $\neg$ CH.

\section{Products of torsion modules}

Our first aim is to prove that every  weakly productively
bounded domain  is actually productively bounded.  We start
by recording a lemma from \cite{gz}  to which we will refer
repeatedly.

\begin{lemma} {\rm (see \cite[Lemma 1.1]{gz})}
\label{cip}If $R$ is weakly productively bounded then, for
every  uncountable family ${\cal A}$ of ideals such that
$\cap {\cal A}=0$, there is a countable subfamily ${\cal A}
^{\prime }$ of ${\cal A}$ such that $\cap {\cal A}^{\prime
}=0$. $\Box $
\end{lemma}

 Next we list some consequences of the assumption that there
is a counterexample to our claim. Call a domain $R$  {\it
ultrafiltral} if there is a non-principal ultrafilter $U$ on
$\omega $,  together with a family
$(A_n)_{n \in  \omega }$ of ideals of $R$ such that, for
every subset $Y$  of $\omega $, the intersection $\bigcap_{n
\in  Y}A_n$ is zero if and only  if $Y \in  U$. It is
essentially proved in \cite{gz} and \cite{bz} that $R$ is
ultrafiltral in case $R$ is  weakly productively bounded
without being productively bounded; we will sketch the proof
and strengthen the conclusion.

\begin{lemma}
\label{ppt}Suppose that $R$ is weakly productively bounded
but not  productively bounded. Then  there is a {\em p-point}
$U$, together with a  family
$(A_n)_{n \in  \omega }$ of ideals of $R$ such that, for
every subset $Y$  of $
\omega $, the intersection $\bigcap_{n \in  Y}A_n$ is zero if
and only if
$Y \in  U$.
\end{lemma}

{\sc proof. }It follows from the assumptions on $R$ and Lemma
\ref{cip} that there is a {\it countable} family $(M_i)_{i
\in  I}$ of torsion modules such that $\prod_{i \in  I}M_i$
is torsion and $\bigcap_{i \in  I}\limfunc{ann}(M_i)=0
$. Let $A_i=\limfunc{ann}(M_i)$. As in \cite[Thm 6.1, proof
of Step II]{gz}, there is a subset $L$ of $I$ such that
$\bigcap_{i \in  L}A_i=0$ and for any pair of disjoint
subsets $J$ and $K$ of $L$ either $\bigcap_{i \in   J}A_i\neq
0
$ or $\bigcap_{i \in  K}A_i\neq 0$. Without loss of
generality we can assume that $I=L=\omega $. If $U$ is
defined to be $\{Y\subseteq \omega :$ $
\bigcap_{i \in  Y}A_i=0\}$, $U$ is a non-principal
ultrafilter on $\omega $. We claim that $U$  is a p-point.

Suppose that $\omega =\coprod_{k \in  \omega }N_k$ such that
$N_k\notin U$  for every $k \in
\omega $. For each $m \in  \omega $, let
\begin{equation}
\label{tm}T_m=\bigoplus \{M_n:n \in  \bigcup_{k\geq
m}N_k\}\hbox{.}
\end{equation} Then, for all $m$, the annihilator
$\limfunc{ann}(T_m)=\bigcap \{A_n:n \in
\bigcup_{k\geq m}N_k\}$ is zero since  $\bigcup_{k\geq m}N_k
\in  U$ (note that the finite  union $
\bigcup_{k<m}N_k$ does not belong to $U$). Consequently,
since $R$ is  wPB,  the product $\prod_{m \in  \omega }T_m$
is not torsion. Let $y = (y(m))_{m  \in  \omega}$ be an
element of
$\prod_{m \in  \omega }T_m$ which is not of finite order,
i.e., $\bigcap_{m \in  \omega }\limfunc{ann}(y(m))=0$.
Moreover, let
$Y$ be the union of the supports of the $y(m)$; more
precisely
$$ Y=\{n \in  \omega :\exists m\hbox{ s.t. }y(m)\hbox{ has a
non-zero projection on }M_n\}\hbox{.}
$$ (Here we refer to the canonical projections associated
with the definition of $T_m$ in (\ref{tm}).) Then $Y$ belongs
to $U$ since $\bigcap_{n \in  Y}A_n\subseteq
\limfunc{ann}(y)=0$. On the other hand, for all $k \in
\omega $, the intersection $Y\cap N_k$ is clearly finite by
construction. $\Box $

\smallskip\

The contradiction we seek now follows from the following
observation.

\begin{lemma}
\label{bs} Suppose that there is a p-point $U$, together with
a family
$(A_n)_{n \in
\omega} $ of ideals of $R$ such that, for every subset $Y$ of
$\omega
$,  the intersection $\bigcap_{n \in  Y}A_n$ is zero if and
only if $Y$ belongs to $U$. Then there is an  uncountable
family ${\cal   B}$ of ideals of $R$ such that $\bigcap {\cal
B}=0$, while $\bigcap {\cal  B}^{\prime}\neq 0$ for every
countable subfamily ${\cal B}^{\prime }$ of ${\cal B}$.
\end{lemma}

{\sc proof.} Let $\bar U=\{S\subseteq \omega :S$ is infinite
and $S\notin U\}
$. For each $S \in  \bar U$ and $m \in  \omega $ we moreover
define
$$ B_{S,m}=\bigcap \{A_n:n \in  S\setminus
\{0,1,...,m\}\}\hbox{.}
$$ Note that the $B_{S,m}$ form an ascending chain of ideals
each of which is non-zero because $S$ (and hence $S\setminus
\{0,1,...,m\}$) does not belong to
$U$.

Let $B_S=\bigcup_{m \in  \omega }B_{S,m}$, and set ${\cal
B}=\{B_S:S \in  \bar U\}$. First we claim that $\bigcap {\cal
B}=0$. Suppose, to the contrary,  that $\bigcap {\cal B}$
contains a nonzero element $r$. If $X=\{n \in
\omega :r \in  A_n\}$, then clearly $X$ does not belong to
$U$ since $r \in  \bigcap_{n \in   X}A_n$. Therefore $
\omega \setminus X \in  U$, and in particular, $\omega
\setminus X$ is  infinite. Whenever we write $\omega
\setminus X$  as the disjoint union of two infinite subsets,
these cannot both belong to $U$.  Hence there is a subset
$S_1\subseteq
\omega \setminus X$ which belongs to $\bar U$. But then $r
\in  B_{S_1}$  and consequently $r \in  B_{S_1,m}$ for some
$m$.  This means that $S_1\setminus \{0,1,...,m\}$ is
contained in $X$, a contradiction to the choice of $S_1$.

It remains to be proved that $\bigcap {\cal B}^{\prime }\neq
0$ for every  countable subset ${\cal B}^{\prime }$ of $
{\cal B}$.  Observe that so far we have only used the
ultrafilter  properties of $U$.  Now we will use the fact
that $U$ is a p-point to show that,  whenever $\{S_n:n \in
\omega \}\subseteq \bar U$,
  there exists $S \in  \bar U$ such that $B_S\subseteq \bigcap
\{B_{S_n}:n \in  \omega \}$. This will clearly imply our
claim concerning  countable subfamilies of ${\cal B}$.

If $\bigcup_{n \in  \omega }S_n \in  \bar U$, we can take
$S=\bigcup_{n \in  \omega }S_n$. So assume $\bigcup_{n \in
\omega }S_n \in  U $. Let $N_0=\omega \setminus \bigcup_{n
\in  \omega }S_n$, $N_1=S_0$, and for $k>1$, let
$N_k=S_{k-1}\setminus \bigcup_{n<k-1}S_n$. Then $N_k\notin
U$ for all $k \in
\omega $, and $\omega$ is the disjoint union of the $N_k$.
Hence, there exists $Y \in  U$ such that, for all $k \in
\omega $, the intersection
$Y\cap N_k$ is finite. Let $S=\omega \setminus Y$. Clearly
$S$ is infinite, since $Y\cap S_0$ is finite, and so $S \in
\bar U$.

To verify the claimed inclusion, we need to show that
$B_{S,m}\subseteq  B_{S_n}$ for all $m,n \in  \omega $. Our
construction permits us to write $S_n=\tilde S_n\cup F_n$ for
$n  \in
\omega$, where $\tilde S_n\subseteq S$ and $F_n\subseteq
Y\cap \bigcup_{k\leq n+1}N_k$ is finite. This clearly implies
$\tilde S_n\setminus \{0,1,...,m\}\subseteq S\setminus
\{0,1,...,m\}$ for all $m \in  \omega $. Therefore
$$
\begin{array}{c}
\bigcap \{A_\ell :\ell  \in  S\setminus
\{0,1,...,m\}\}\subseteq \bigcap
\{A_\ell :\ell  \in  \tilde S_n\setminus \{0,1,...,m\}\} \\
=\bigcap \{A_\ell :\ell  \in  S_n\setminus (F_n\cup
\{0,1,...,m\})\}\subseteq B_{S_n},
\end{array}
$$ which completes the proof of the Lemma. $\Box $

\smallskip\

\begin{theorem}
\label{pb} Every weakly productively bounded domain is
productively  bounded. \end{theorem}

{\sc proof. }Suppose to the contrary that $R$ is wPB but not
PB. By Lemma
\ref{ppt} the hypothesis of  Lemma \ref{bs} is satisfied. But
then $R$ is  not wPB by Lemma \ref{cip}, a contradiction.
$\Box $

\section{Lattices}

Primarily, this section is devoted to proving that, for every
domain $R$,  transversal boundedness of the dual ideal
lattice $L_R$ entails uniform  transversal boundedness of
that lattice.  This does not follow from  Theorem \ref{pb},
because it is unresolved whether every productively  bounded
domain $R$ gives rise to a uniformly transversally bounded
lattice $L_R$; but it does follow from Lemma \ref{bs} and
previous  results of the second author.

As we mentioned earlier, this implication cannot be extended
to general  lattices in the presence of the continuum
hypothesis.  We will  subsequently explore what happens when
the continuum hypothesis fails.

\begin{theorem}
\label{utb}For every domain $R$, transversal boundedness of
$L_R$ is  equivalent to uniform transversal boundedness of
$L_R$.
\end{theorem}

We need to review some of the terminology and results of
\cite{z2} before  being able to apply the insights of the
previous section.  Given any  ultrafilter $U$ on $\omega$, we
consider an associated complete lattice
$L(U) = \{ A \subseteq \omega : A \notin U \} \cup \{
\omega\}$, ordered  by inclusion, in which meets coincide
with set-theoretic intersections  and joins are given by

$$
\bigvee \{A_i:i \in  I\}=\left\{
\begin{array}{ll}
\bigcup_{i \in  I}A_i & \hbox{if }\bigcup_{i \in  I}A_i\notin
U \\ \omega  &
\hbox{otherwise}
\end{array}
\right.
$$

\noindent The following lemma combines two results of
\cite{z2}.

\begin{lemma}
\label{aux} Let $L$ be a complete, finitely join-irreducible
lattice  which is trans\-versally bounded, but not uniformly
so.  Then there exists  a p-point $U$, together with a
complete upper subsemilattice $L'$ of $L$  which (as a
complete upper semilattice) is isomorphic to $L(U)$.

In case $L = L_R$ for a domain $R$, there is a family
$(X_n)_{n  \in  \omega}$ of subsets of $L$ such that
 a semilattice $L'$ as above and an  isomorphism $\phi: L(U)
\rightarrow L'$ can be constructed  as follows:   if $A_n =
\bigcap X_n $, then $L'$ consists of all intersections of
subfamilies of the family $(A_n)_{n  \in  \omega}$, and $\phi
(\{n\}) = A_n$.
\end{lemma}

{\sc proof.}  By Corollary C of \cite{z2}, there exists a
complete upper  subsemilattice $L'$ of $L$ which is
isomorphic to an upper semilattice of  the form $L(U)$ for
some ultrafilter $U$ on $\omega$.  The lattice $L'$,  being
closed under suprema in $L$, clearly inherits the property of
being  TB, and hence $L(U)$ is TB.  But, by Theorem E(I) of
\cite{z2}, this  guarantees that $U$ is a p-point.

The claim concerning a realization of $L'$ in case $L$ is the
dual ideal  lattice of a domain, is an immediate consequence
of part (c) of Corollary  C of \cite{z2}.  $\Box$

\smallskip\

{\sc proof of Theorem \ref{utb}.} Suppose, to the contrary,
that $L_R$ is TB  without being UTB, and let $U$, $L'$,
$(X_n)_{n  \in  \omega}$ and $A_n$ be  as in Lemma
\ref{aux}.  In particular, $U$ is a p-point.  Moreover, the
definition of $L(U)$ immediately yields the following string
of  equivalences for any $Y \subseteq \omega$:   $Y  \in  U$
if and only  if \, $\bigvee Y = \omega$ in $L(U)$ if and only
if \, $\phi(\bigvee Y) = 0$.   But $\phi(\bigvee Y) =
\bigcap_{n  \in  Y} A_n$, and so the p-point $U$ and  the
family $(A_n)_{n  \in  \omega}$ of ideals satisfy  the
hypothesis of  Lemma \ref{bs}.  The conclusion of Lemma
\ref{bs}, when combined with  Lemma \ref{cip},  shows that
$R$ fails to be wPB.  On the other hand,  transversal
boundedness of $L_R$ clearly forces $R$ to be wPB. This
contradiction completes the proof. $\Box$

\smallskip\

The second author  has shown that CH implies the existence of
a complete finitely join-irreducible lattice which is
transversally bounded  without having the uniform boundedness
property (\cite[p. 204]{z2}). Here we observe that this
conclusion is independent of $\neg $CH.

\begin{theorem}
\label{meta}It is undecidable in ZFC + $\neg $CH whether
there is a complete finitely join-irreducible lattice which
is transversally bounded, but not  uniformly so.
\end{theorem}

{\sc proof.}  The axioms MA + $\neg $CH imply that there is a
Ramsey ultrafilter (see for example \cite[p.259]{J}). Hence,
by Example J and Theorem G of \cite{z2}, any model of ZFC +
MA + $\neg $CH  admits a complete finitely join-irreducible
lattice which is TB but not UTB.

On the other hand, in a model of ZFC + $\neg $CH without
p-points, there  is no such lattice by Lemma \ref{aux}.
$\Box $

\section{Domains with Krull dimension}

For domains with Krull dimension we can show the equivalence
of the ring-theoretic and the lattice-theoretic properties
considered in the previous sections. Moreover, in this case,
all of these boundedness  conditions follow from the (on the
face of it comparatively weak)  condition on uncountable
families of ideals which arises as a consequence of weak
productive  boundedness in Lemma \ref{cip}.  More precisely,
we have:

\begin{theorem}
\label{krull}Suppose $R$ is a domain with Krull dimension and
$L_R$ its  dual ideal lattice. Then the following statements
are equivalent:

(1) $R$ is weakly productively bounded;

(2) $R$ is productively bounded;

(3) $L_R$ is transversally bounded;

(4) $L_R$ is uniformly transversally bounded;

(5) for every uncountable family ${\cal A}$ of ideals of $R$
such that
$\cap  {\cal A}=0$, there is a countable subfamily ${\cal
A}^{\prime }$ of
${\cal A}
$ such that $\cap {\cal A}^{\prime }=0$.
\end{theorem}

\noindent {\it Remark.}  Our argument for the crucial
implication `(5)
$\Longrightarrow $ (4)' was inspired by the proof of Theorem
8 in
\cite{bg}.  We recall the pivotal definition introduced
there.  Given a  complete lattice $L$, we start by fixing a
family $(X_i)_{i  \in  I}$ of  subsets of $L$.  An element
$y  \in  L$ is said to be {\it tame} relative to an  element
$x  \in  L$ if there exists an infinite subset $K \subseteq
I$,  together with a family of elements $x_k  \in  X_k$ for
$k  \in  K$, such that
$x \vee \bigvee_{k  \in  S} x_k \ge y$ for every infinite
subset $S$ of $K$. For the convenience of the reader, we
include Lemma 7 of \cite{bg}.

\begin{lemma}
\label{bergg}  Let $L$ and $(X_i)_{i  \in  I}$ be as above
and suppose that
$L$ does not contain a complete upper subsemilattice
isomorphic to
$2^{\omega}$.  Then there exists a finite subset $F \subseteq
I$,  together with a family $(x_i)_{ i  \in  F}$ of elements
$x_i  \in  X_i$, such  that every element of \ $\bigcup_{i
\in  I \setminus F} X_i$ is tame  relative to $\bigvee_{i
\in  F} x_i$.  $\Box$ \end{lemma}

{\sc proof of Theorem \ref{krull}. } The implications `(4)
$\Longrightarrow $ (2)' and `(3) $
\Longrightarrow $ (1)' are known, and `(2) $\Longrightarrow $
(1)', as  well as `(4) $
\Longrightarrow $ (3)' are trivial.  Moreover, `(1)
$\Longrightarrow $  (5)' follows from Lemma \ref{cip}. Hence
it suffices to prove `(5) $\Longrightarrow $ (4)'.  We assume
(5) and let
$(X_i)_{i  \in  I}$ be a family of nonempty subsets of the
dual ideal  lattice $L_R$ such that, for each cofinite subset
$J \subseteq I$, the  intersection of the ideals in
$\bigcup_{i \in  J} X_i$ is zero.  We wish to  apply Lemma
\ref{bergg} to construct a transversal of $(X_i)_{i  \in
I}$  which is unbounded in $L_R$.  Clearly, it is harmless to
assume that none  of the sets $X_i$ contains the zero ideal.

Since $R$ has Krull dimension, $L_R$ does not contain any
subsets  order-isomorphic to $2^{\omega}$, and consequently
$L_R$ satisfies the  hypothesis of Lemma \ref{bergg}.  We
infer the existence of a finite  subset $F \subseteq I$ and a
family of ideals $(A_i)_{i  \in  F}  \in
\prod_{i  \in  F} X_i$ such that each ideal $A  \in
\bigcup_{i  \in  I
\setminus F} X_i$ is tame relative to the intersection $B =
\bigcap_{i
 \in  I} A_i$, inside the dual ideal lattice of $R$.

Our initial assumption about the $X_i$ forces the
intersection of the  ideals in $\bigcup_{i  \in  I \setminus
F} X_i$ to be zero, and hence  Condition (5) provides us with
a countable family $(A_n)_{n  \in  \omega}$  of ideals in
$\bigcup_{i  \in  I \setminus F} X_i$ such that $\bigcap_{n
 \in  \omega} A_n = 0$.  Since each $A_n$ is tame relative to
$B$, we can  moreover find infinite subsets $K_n \subseteq I$
and transversals
$(A_{k,n})_{k  \in  K_n}  \in  \prod_{k  \in  K_n} X_k$ such
that $$ B\cap
\bigcap_{k  \in  S_n} A_{k,n} \subseteq A_n$$ for each
infinite subset $S_n \subseteq K_n$.  A standard diagonal
technique then yields a family $(S_n)_{ n  \in  \omega}$ of
pairwise  disjoint infinite sets $S_n \subseteq I$ such that
$S_n$ is contained in
$K_n$ for $n  \in  \omega$. Moreover, we define $S_{-1} = I
\setminus
\bigcup_{n
 \in
\omega} S_n$ and, for each $i  \in  S_{-1}$, we pick an
arbitrary ideal $A_{i, -1}
 \in  X_i$.  We will check that the transversal $(A_{k,n})_{k
 \in  S_n,  n \ge -1}$ of \, $\prod_{i  \in  I} X_i$ is
unbounded in $L_R$. Indeed, our construction entails that $$B
\cap \bigcap_{k  \in  S_n, n \ge -1} A_{k, n} \subseteq B
\cap \bigcap_{k  \in  S_n, n \ge 0}  A_{k, n}
\subseteq \bigcap_{n  \in  \omega} A_n = 0.$$ But since $B$,
being a finite intersection of nonzero ideals, is nonzero,
this implies that
$\bigcap_{k  \in  S_n, n \ge -1}  A_{k, n} = 0$ as desired.
Thus  $L_R$ is UTB, that is, condition (4) is satisfied.
$\Box$

\section{Countable extensions of domains}

Our aim is to prove the following, which generalizes (for the
commutative case) Theorems 4.2 and 6.1 of \cite{gz}.\

\begin{theorem}
\label{exts2} Suppose that $D\subseteq R$ is an extension of
commutative  domains such that $R$ is countably generated
over $D$.  If $D$ is  productively bounded then so is $R$.
\end{theorem}

\noindent {\it Remark.}  The referee has pointed out that the implication of
this theorem remains valid if $D$ is a right Ore domain and
the extension $R$ is generated over $D$ by countably many
elements which are central in $R$.
\medskip

{\sc proof. } Let $D$ be PB.  By Theorem \ref{pb}, it is
enough to prove  that  $R$ is wPB, i.e.,  to show that if $
(M_i)_{i \in  I }$ is a family of $R$-modules with
$\limfunc{ann}_R(M_i)=0
$ for all $i \in  I $, then $\prod_{i \in  I }M_i$ is not
torsion. Clearly,  without loss of generality, we may assume
that our family is countable,  so $I$ may be taken to be
$\omega$.

We first observe that we can reduce the situation to the case
where $R$  is finitely generated over $D$ as a ring, as
follows. If $R=D[x_n:n \in  \omega ]$,  let
$R_k=D[x_n:n<k]$.  Moreover, we write $\omega = \coprod_{k
\in
\omega} N_k$, where the $N_k$ are pairwise disjoint infinite
sets.   Assuming the result for the finitely generated case,
we obtain, for each $k$, an  element $(m_n)_{n \in  N_k}$ in
$\prod_{n \in  N_k}M_n$ which is not torsion  over $R_k$.
Clearly the transversal $(m_n)_{n \in  \omega }  \in
\prod_{n \in  \omega }M_n$ then fails  to be a torsion element
over $R$. (Note that the proof shows that the union of a
countable  chain of productively bounded domains is
productively bounded.)

By induction the finitely generated case can obviously be
reduced to the  following: if
$D$ is PB and $R=D[x]$, then $R$ is PB. So suppose that
$R=D[x]$.  We say that an element of $R$ has degree $\leq k$
iff it can be written in the  form
$\sum_{i=0}^kd_ix^i$ for some $d_i \in  D$. Once more, we
decompose
$\omega$ into infinitely many disjoint infinite subsets
$N_k$.  It clearly  suffices to prove that, for every $k \in
\omega $, there is an element $(m_n)_{n \in  N_k}$ in
$\prod_{n \in  N_k}M_n$ such that $\bigcap
\{\limfunc{ann}_R(m_n):n \in  \bigcup_{\ell \leq k}N_\ell \}$
does not contain any non-zero elements of degree $\leq k$.
We construct such transversals $(m_n)_{n \in  N_k}$ by
induction on $k$.   For $k=0$, we obtain $  (m_n)_{n \in
N_0}$ as required from the fact that $D$ is wPB;  indeed,
the  latter implies that $\prod_{n \in  N_0}M_n$ is not
torsion as a $D$-module.

Suppose now that elements $(m_n)_{n \in  N_\ell }$ with the
desired  properties have been defined for $\ell \leq k$. Set
$A_k=\bigcap
\{\limfunc{ann}_R(m_n):n \in
\bigcup_{\ell \leq k}N_\ell \}$.  In particular, $A_k$ then
does not  contain any nonzero elements of degree $\leq k$.
If $A_k$ does not  contain any elements of degree $\leq k+1$,
let $(m_n)_{n \in  N_{k+1}}$ be the zero element. Otherwise,
pick an element
$$ r_{k+1}=d_{k+1}x^{k+1}+t_k
$$ in $A_k$, where $d_{k+1} \in  D\setminus \{0\}$ and $t_k
\in  R$ has degree $
\leq k$. Since $\limfunc{ann}_D(r_{k+1}M_n)=0$ for all $n
\in  N_{k+1}$ and $D$ is wPB, there is an element
$(r_{k+1}m_n)_{n \in  N_{k+1}} \in  \prod_{n \in
N_{k+1}}M_n$ which is not $D$-torsion. Hence
\begin{equation}
\label{eqn}Dr_{k+1}\cap \bigcap \{\limfunc{ann}_R(m_n):n \in
N_{k+1}\}=0\hbox{.}
\end{equation} It follows that $\bigcap
\{\limfunc{ann}_R(m_n):n \in  \bigcup_{\ell \leq k+1}N_\ell
\}=A_k\cap \bigcap \{\limfunc{ann}_R(m_n):n \in  N_{k+1}\}$
does not contain any elements of degree $\leq k+1$. Indeed,
if $r=d x^{k+1}+t$ were such an element in the intersection
(where $d  \in  D$ and
$t$ has degree $\leq k$), then $  d_{k+1}r=d
d_{k+1}x^{k+1}+dt_k$ because there is at most one element of
the form $d d_{k+1}x^{k+1}+t^{\prime }$ in $A_k$ with $
t^{\prime }$ of degree $\leq k$;  keep in mind that $A_k$
contains no  non-zero elements of degree $\leq k$. Thus
$d_{k+1}r \in  Dr_{k+1}$, which contradicts (\ref{eqn}).

This completes the inductive step, and hence the proof. $\Box
$

\smallskip

Combining Theorem \ref{exts2} with a result of Bergman and
Galvin, we obtain

\begin{corollary} Every commutative domain which is countably
generated over  a noetherian subring is productively bounded.
\end{corollary}

{\sc proof.} By Theorem \ref{exts2} above and Corollary 11 of
\cite{bg}. $\Box $

\smallskip\

In a similar manner one can prove:

\begin{theorem}
\label{utb2} Suppose $D\subseteq R$ is an extension of
commutative  domains such that $R_D$ is countably generated.
If the dual ideal lattice of $D$ is uniformly  transversally
bounded, then the same is true for the dual ideal lattice  of
$R$. $\Box $
\end{theorem}

As a consequence of Theorems \ref{krull} and \ref{utb2}, we
have:

\begin{corollary} Suppose $D\subseteq R$ are commutative
domains such that $R$ is countably  generated over $D$ and
$D$ has Krull dimension. Then weak productive  boundedness of
$D$ implies uniform transversal boundedness of the lattice
$L_R$. $\Box $
\end{corollary}

\bigskip\ \

Note: Prior to 1995, the second author published under the
name `Zimmermann-Huisgen'.

%end of body

\end{document}